\documentclass[final]{siamart171218}

\usepackage{lipsum}
\usepackage{amsfonts}
\usepackage{graphicx}
\usepackage{epstopdf}
\usepackage{algorithmic}
\ifpdf%
  \DeclareGraphicsExtensions{.eps,.pdf,.png,.jpg}
\else
  \DeclareGraphicsExtensions{.eps}
\fi
\usepackage{amsopn}
\DeclareMathOperator{\diag}{diag}
\usepackage{booktabs}

\usepackage{mathtools}
\usepackage{bm}
\usepackage{amssymb}

\DeclareMathOperator*{\argmin}{arg\,min}

\newcommand{\TheTitle}{%
	On selecting coarse-grid operators for Parareal and MGRIT applied to linear advection
}

\newcommand{\TheShortTitle}{%
	Parareal and MGRIT convergence for the advection equation
}

\newcommand{\TheName}{%
  O. A. Krzysik
}

\newcommand{\TheAddress}{%
  School of Mathematical Sciences, Monash University, Clayton, Australia
  (\email{oliver.krzysik@monash.edu}).
}

\newcommand{\TheFunding}{%
	This work was submitted to the Student Paper Competition of the 19th Copper Mountain Conference On Multigrid Methods, Copper Mountain, CO, USA. 
	
	Much of this manuscript has formed the foundations of a full-length manuscript for which a pre-print is available at \url{https://arxiv.org/abs/1910.03726}.
	
  	This work was funded by an Australian
Government Research Training Program (RTP) Scholarship.
}

\newcommand{\TheCollaborators}{%
  H. De Sterck, S. P. MacLachlan, S. Friedhoff
}

\author{\TheName\thanks{\TheAddress}}
\title{{\TheTitle}\thanks{\TheFunding}}
\headers{\TheShortTitle}{\TheName}
\ifpdf%
\hypersetup{%
  pdftitle={\TheTitle},
  pdfauthor={\TheName}
}

\begin{document}

\maketitle

\begin{center}
In collaboration with:
  {\TheCollaborators}
\end{center}
\vspace{1cm}

\begin{abstract}
We consider the parallel time integration of the linear advection equation with the Parareal and two-level multigrid-reduction-in-time (MGRIT) algorithms. Our aim is to develop a better understanding of the convergence behaviour of these algorithms for this problem, which is known to be poor relative to the diffusion equation, its model parabolic counterpart. Using Fourier analysis, we derive new convergence estimates for these algorithms which, in conjunction with existing convergence theory, provide insight into the origins of this poor performance. We then use this theory to explore improved coarse-grid time-stepping operators. For several high-order discretizations of the advection equation, we demonstrate that there exist non-standard coarse-grid time stepping operators that yield significant improvements over the standard choice of rediscretization.
\end{abstract}

\begin{keywords}
	Multigrid, Parareal, MGRIT, parallel-in-time, Fourier analysis, hyperbolic
\end{keywords}

\section{Introduction}\label{sec:intro}
In the context of the large-scale numerical simulation of time-dependent partial differential equations (PDEs), the advent of massively parallel computers, in combination with a desire for faster compute times, has driven the development of highly concurrent algorithms. Relative to traditional methods that use sequential time stepping, these algorithms may make better use of available parallel resources through the use of so-called parallel time integration techniques. Two popular algorithms in this area are Parareal \cite{Lions2001} and multigrid-reduction-in-time (MGRIT) \cite{Falgout2014}. Parallel testing for the diffusion equation has shown these algorithms can provide significant speed-ups (in wall clock time) over sequential time marching, given enough parallel resources \cite{Falgout2014, Falgout2017}. 

Despite their success for model parabolic problems, Parareal and MGRIT perform significantly worse for hyperbolic PDEs, or at least PDEs without significant diffusivity \cite{DeSterck2018_mgrit_hyp, Dobrev2017, HowseThesis2017, Nielsen2018, Ruprecht2018, DeSterck2018_analysis}. For example, in \cite{DeSterck2018_mgrit_hyp} relatively little parallel speed-up was attained for linear advection with simple 1st-order discretizations despite the use of significant parallel resources. The focus of the current work is developing a better understanding of this poor convergence behaviour for the linear advection equation. There has already been significant progress made towards understanding the general convergence properties of these algorithms. In \cite{Friedhoff2015}, a semi algebraic mode analysis (SAMA) methodology was proposed for estimating convergence rates of multigrid methods, and applications for MGRIT were shown. Convergence bounds were derived for two-level MGRIT in \cite{Dobrev2017}, and under quite general assumptions, convergence of Parareal and multilevel MGRIT was proven in \cite{Southworth2018}. Here, we derive new convergence estimates using Fourier analysis, and develop optimization strategies for the selection of coarse-grid operators, yielding significantly improved convergence rates.

\section{Preliminaries} \label{sec:prelims}
Here, the model problem is outlined, and a brief overview of Parareal and two-level MGRIT is given.

\subsection{Model problem}
\label{ssec:model_problem}

We are concerned with the numerical solution of the one-dimensional linear advection equation,
\begin{align} \label{eq:LA_PDE}
\partial_t u + \partial_x (au) = 0, 
\quad 
u(x,0) = g(x), 
\quad (x,t) \in \Omega \coloneqq (-1,1) \times (0, T_{\rm f}),
\end{align}
subject to periodic spatial boundary conditions, $u(-1,t) = u(1,t)$, and with constant wavespeed $a > 0$. Despite \eqref{eq:LA_PDE} being the simplest example of a hyperbolic PDE, it has proven difficult to solve efficiently using Parareal and MGRIT, even for simple 1st-order discretizations, \cite{DeSterck2018_mgrit_hyp, HowseThesis2017, Ruprecht2018}, which provides the motivation for its consideration here.

We define a uniform spatial mesh $\bm{x} \coloneqq (x_j)_{j = 0}^{n_x} = (-1 + j \Delta x)_{j = 0}^{n_x}$, with $\Delta x = 2/n_x$, and a uniform temporal mesh $\bm{t} \coloneqq (t^n)_{n=0}^{n_t} = (n \Delta t)_{n=0}^{n_t}$, with $\Delta t = T_{\rm f}/ n_t$. Equation \eqref{eq:LA_PDE} is discretized in space, and the resulting ordinary differential equations are discretized in time using a one-step method. Letting $\bm{u}^n \in \mathbb{R}^{n_x}$ denote the discrete spatial solution at time $t^n$, we write the fully discretized problem as
\begin{align} \label{eq:discrete_PDE_general}
\bm{u}^{n+1} = \Phi \bm{u}^n,  \quad \bm{u}^0 = g(\bm{x}),  
\quad n = 0,\dots, n_t-1,
\end{align}
with $\Phi \in \mathbb{R}^{n_x \times n_x}$ known as the time-stepping operator, or simply as the time stepper.

In particular, for results shown in \S \ref{sec:psi_choice}, we use the 2nd- and 3rd-order upwind-finite-difference spatial discretizations given by
\begin{align}
\big[
\partial_x
(au)
\big]_j = 
\begin{cases}
\displaystyle{\frac{a}{2 \Delta x}
\big[3u_{j} - 4u_{j-1} + u_{j-2} \big] 
+ 
\mathcal{O}(\Delta x^2)}, 
\\[10pt]
\displaystyle{\frac{a}{6 \Delta x}
\big[
u_{j-2} - 6u_{j-1} + 3u_j + 2u_{j+1}
\big]
+
 \mathcal{O} \left( \Delta x^3 \right)}.
\end{cases}
\end{align}
Temporal discretizations are carried out with three Runge--Kutta (RK) schemes: 2nd-order Heun's method; the optimal, 3rd-order strong-stability preserving method; and a 3rd-order L-stable SDIRK method (see \cite[p. 262]{Butcher2003}). These RK schemes are, respectively, defined by the following Butcher tableaus:
\renewcommand*{\arraystretch}{1.3}
\begin{align} \label{eq:rk_butcher}
\begin{array}
{c|ccccc}
0 & 0 & 0 & \\
1 & 1 & 0 & \\
\hline
& \tfrac{1}{2} & \tfrac{1}{2}
\end{array}
\quad 
\quad
\begin{array}
{c|ccccc}
0   & 0     & 0    & 0  \\
1    & 1     & 0    & 0  \\
\tfrac{1}{2} & \tfrac{1}{4} & \tfrac{1}{4} & 0  \\
\hline
& \tfrac{1}{6} & \tfrac{1}{6} & \tfrac{2}{3}
\end{array}
\quad 
\quad
\begin{array}
{c|ccccc}
\zeta  & \zeta   & 0  & 0  \\
\alpha   & \beta  & \zeta  & 0  \\
1   & \gamma  & \epsilon & \zeta  \\
\hline
& \gamma & \epsilon & \zeta
\end{array}
\end{align}
with $\zeta = 0.4358665215\dots$, $\alpha = \tfrac{1+\zeta}{2}$, $\beta = \tfrac{1-\zeta}{2}$, $\gamma = -\tfrac{3}{2}\zeta^2 + 4\zeta - \tfrac{1}{4}$, and $\epsilon = \tfrac{3}{2}\zeta^2 - 5\zeta + \tfrac{5}{4}$. The 2nd-order spatial discretization is paired with the 2nd-order RK method to create an explicit 2nd-order scheme with CFL condition $a \tfrac{\Delta t}{ \Delta x} \leq \tfrac{1}{2}$. The 3rd-order spatial discretization is combined with the 3rd-order RK methods to form 3rd-order explicit and implicit schemes. The 3rd-order explicit scheme has the CFL condition $a \tfrac{\Delta t}{ \Delta x} \lesssim 1.625$, while the 3rd-order implicit scheme is unconditionally stable. For these temporal discretizations, it can be shown that each time stepper $\Phi$ is a rational function in its respective spatial discretization. Additionally, since periodic spatial boundary conditions are used, $\Phi$ will be circulant.

\subsection{Parareal and MGRIT}
MGRIT \cite{Falgout2014} is a multigrid algorithm for solving equations in the form of \eqref{eq:discrete_PDE_general} (which constitute a block lower triangular linear system) through approximating a block cyclic reduction procedure. Here we give a brief overview of the two-level MGRIT variant, which is also related to Parareal \cite{Lions2001}; for more detailed descriptions, see \cite{Dobrev2017, Falgout2014}.

We first define a coarse temporal mesh through an integer coarsening factor $m > 1$: $(T^n)_{n = 0}^{n_t/m} \equiv (n \Delta T)_{n = 0}^{n_t/m}$, with $\Delta T = m \Delta t$. The coarse time points are referred to as C-points, and the set of points not appearing on the coarse mesh are F-points.  

Like any typical multigrid algorithm, MGRIT combines relaxation with a coarse-grid correction procedure. The two fundamental types of relaxation are: F-relaxation, which is time stepping from each C-point across the following F-interval; and C-relaxation, which is time stepping from the last F-point in each interval to its following C-point. The standard relaxation sweeps performed in MGRIT are either: F-relaxation (for which the algorithm is equivalent to a Parareal-type algorithm), or FCF-relaxation, which is an F-, followed by a C-, followed by an F-relaxation. The interpolation operator is motivated through the reduction framework, and is known as ideal interpolation. It is equivalent to injection interpolation at C-points followed by an F-relaxation. The restriction operator is simply injection restriction.

For the coarse-grid correction problem, the error, $\bm{e}$, at C-points, $t^n = m n \Delta t$, $n = 0, \dots n_t/m$, is approximated by the system
\begin{align} \label{eq:cg_error_equation}
\bm{e}^{m(n+1)} = \Psi \bm{e}^{mn} + \bm{r}^{m(n+1)},  \quad \bm{e}^0 = 0, \quad n = 0,\dots, n_t/m  -1,  \quad \bm{e}^{n} \in \mathbb{R}^{n_x},
\end{align}
with $\bm{r}^{m(n+1)}$ the fine-grid residual at the $n+1$st C-point. Here, $\Psi \in \mathbb{R}^{n_x \times n_x}$ is the coarse-grid time stepper. Taking $\Psi = \Phi^m$ results in a Schur complement coarse-grid operator, yielding an exact two-level algorithm; however, the resulting coarse-grid problem is as expensive to solve as the fine-grid one, and so it is not considered practically feasible. Instead, one chooses $\Psi \approx \Phi^m$; often this approximation is made through rediscretization of $\Phi$ on the coarse grid. However, what makes a good choice for $\Psi$ is still very much an open question, and is the subject of \S \ref{sec:convergence} and \S \ref{sec:psi_choice}. 

%
%

\section{Convergence theory} \label{sec:convergence}
With the aim of better understanding the convergence of MGRIT, we now analyse the error propagator, or iteration matrix, $T$, of the algorithm. The error propagation matrix describes the evolution of an initial error under the action of the algorithm. That is, given an initial space-time error $\bm{e}^{(0)}$, after $k$ MGRIT iterations, the error obeys $\Vert \bm{e}^{(k)} \Vert = \Vert T^{k} \bm{e}^{(0)} \Vert \leq \Vert T^{k} \Vert \Vert \bm{e}^{(0)} \Vert  \leq \Vert T \Vert^k \Vert \bm{e}^{(0)} \Vert$. Due to the nilpotency of $T$, the short term convergence behaviour of the algorithm is much more accurately reflected by its norm rather than its spectral radius, which is zero. In this section, we use Fourier analysis \cite{Trottenberg2001} to approximate $\Vert T \Vert_2$, and we also consider the error estimates of \cite{Dobrev2017}. Following this, we give a general discussion on what is required of the coarse-grid time stepper $\Psi$ to achieve fast MGRIT convergence.

\subsection{Fourier analysis}
In numerical experiments (not shown here), we observe for advection-reaction PDEs (corresponding to adding a term to \eqref{eq:LA_PDE} that is proportional to $u$) that short-term MGRIT convergence tends to be remarkably similar for both time-periodic and initial-value problems. We use this as motivation for analysing $T$ with Fourier analysis in time, which applies rigorously (i.e., exactly) in the time-periodic setting, and provides an asymptotically accurate approximation to the initial-value setting considered here (as $n_t \to \infty$). Due to these numerical results, we hypothesize that the Fourier analysis approximation of the initial-value problem is relatively tight. Here, to derive new estimates of $\Vert T \Vert_2$, we apply rigorous Fourier analysis in space, and then local Fourier analysis (LFA) in time.

In what follows, $A$ denotes the fine-grid operator involving $\Phi$ that arises from writing \eqref{eq:discrete_PDE_general} in block matrix form, $A_{\rm c}$ denotes its analogue on the coarse grid involving $\Psi$ that arises from writing \eqref{eq:cg_error_equation} in block matrix form, and $P_{\Phi}$ and $R_I$ are ideal interpolation and injection restriction operators, respectively. The two-grid iteration matrix is \cite{Friedhoff2015, DeSterck2018_analysis}
\begin{align} \label{eq:MGRIT_error_prop}
T 
= 
\left[ I - P_{\Phi} (A_{\textrm{c}})^{-1} R_I A \right] S 
\in 
\mathbb{R}^{n_x n_t \times n_x n_t},
\end{align}
with $S$ the iteration matrix of the smoother, whether it be F- or FCF-relaxation.

Rigorous Fourier analysis \cite{Trottenberg2001} is now applied to the spatial components of $T$. Whilst $\Phi$ is a rational function of the spatial discretization, and so is diagonalizable by its eigenvectors, we place no restriction on $\Psi$ being a function of the fine-grid spatial discretization. However, we do require that $\Psi$ is circulant (for the sake of enforcing periodic spatial boundary conditions on the coarse grid), and so fine time stepper $\Phi$ and coarse time stepper $\Psi$ are simultaneously diagonalizable. Due to their circulant structure, $\Phi$ and $\Psi$ are diagonalized by the periodic Fourier modes
\begin{align}
f_{jk} \coloneqq \exp( \imath j \theta_{k} ), 
\quad 
\theta_{k} = \frac{2 \pi k}{n_x}, 
\quad 
k = -\frac{n_x}{2}, \dots,  \frac{n_x}{2} -1,
\end{align}
by way of the unitary eigenvector matrix ${\cal F} \in \mathbb{C}^{n_x \times n_x}$, with $({\cal F})_{jk} = f_{jk}/\sqrt{n_x}$. For future reference, we denote the eigenvalues of $\Phi$ and $\Psi$ as $(\lambda_k)_{k = -n_x/2}^{n_x/2-1}$ and $(\mu_k)_{k = -n_x/2}^{n_x/2-1}$, respectively.
Following \cite{Friedhoff2015}, we use ${\cal F}$ to diagonalize the $\Phi$ and $\Psi$ blocks present in \eqref{eq:MGRIT_error_prop} and then permute the resulting matrix from spaceline to timeline ordering. This results in a similar block diagonal error propagator
\begin{align}
\widetilde{T} 
\coloneqq 
{\cal P}^\top \left(I_{n_t} \otimes {\cal F}^* \right) T \left( I_{n_t} \otimes {\cal F} \right) {\cal P}  
=
\diag 
\left( 
\widetilde{T}_{-n_x/2}, \dots , \widetilde{T}_{n_x/2-1} 
\right)
\in 
\mathbb{C}^{n_x n_t \times n_x n_t},
\end{align}
with Toeplitz diagonal blocks given by
\begin{align} \label{eq:MGRIT_T_tilde_block}
\widetilde{T}_k
=
\left[ 
I_{n_t} -  \widetilde{P}_{\Phi,k} \left( \widetilde{A}_{{\rm c}, k} \right)^{-1}  \widetilde{R}_{I} \widetilde{A}_k 
\right] \widetilde{S}_k
\in 
\mathbb{C}^{n_t \times n_t}.
\end{align}
Here, $\widetilde{A}_k$ and $\widetilde{A}_{{\rm c}, k}$ are bidiagonal Toeplitz matrices with unit diagonal, and $-\lambda_k$ and $-\mu_k$ on their subdiagonals, respectively. 

We now apply LFA in time to approximately block diagonalize the Toeplitz blocks \eqref{eq:MGRIT_T_tilde_block}, which are not diagonalizable by periodic Fourier modes since they are not circulant (as they are in the time-periodic setting). Before proceeding, for simplicity we limit ourselves to consider only factor-two coarsening ($m = 2$); arbitrary coarsening factors $m$ may be considered, as in \cite{DeSterck2018_analysis}, for example. The problem is extended to the infinite temporal grid $\bm{G}_{\Delta t} \coloneqq \left\{t = \ell \Delta t \, \mid \, \ell \in \mathbb{N}_0 \right\}$, on which we define the continuous coarse-fine (CF) Fourier modes
\begin{align} \label{eq:LFA_CF_modes}
g_{\vartheta, t}^{\rm C} \coloneqq \exp \left(
\imath \vartheta t / \Delta t
\right), \quad 
t \in \bm{G}_{2 \Delta t}; 
\quad 
g_{\vartheta, t}^{\rm F} \coloneqq \exp \left(
\imath \vartheta t / \Delta t
\right), 
\quad
t \in \bm{G}_{\Delta t} \setminus \bm{G}_{2 \Delta t},
\end{align}
with frequency $\vartheta \in \left[ - \pi, \pi \right)$ varying continuously. 

These modes are equivalent to standard red--black Fourier modes, where red points are associated with C-points, and black points with F-points. For each $k$, we extend \eqref{eq:MGRIT_T_tilde_block} to an infinite Toeplitz matrix, and then permute it from natural timeline ordering to CF timeline ordering, in which C-points are blocked before F-points. We then diagonalize the resulting matrix using the continuous Fourier modes \eqref{eq:LFA_CF_modes}. Permuting back to the original ordering, we obtain an infinite, block diagonal matrix with the diagonal blocks:
\begin{align} 
\widehat{T}^{\textrm{F}, \vartheta}_k 
&=
\frac{\lambda^2_k - \mu_k}{\exp(\imath \vartheta) - \mu_k}
\begin{bmatrix}
		1			& 0 \\
\lambda_k 	& 0
\end{bmatrix} 
\in \mathbb{C}^{2 \times 2}
\quad 
&\textrm{(for F-relaxation),}
\label{eq:TPP_T_2x2_block_F}
\\ 
\widehat{T}^{\textrm{FCF}, \vartheta}_k 
&= 
\lambda^2_k \exp({-\imath \vartheta})\, \widehat{T}^{\textrm{F}, \vartheta}_k 
\in \mathbb{C}^{2 \times 2}
\quad 
&\textrm{(for FCF-relaxation)}.
\label{eq:TPP_T_2x2_block_FCF}
\end{align}
The singular values of these $2 \times 2$ matrices are easily calculated, and their maximum over $\vartheta$ gives the two-norm of the respective infinite matrix, which is an approximation to the two-norm of its finite-dimensional counterpart. In each case, one singular value is zero while the other can be expressed and bounded as
\begin{align}
\left\Vert
\widehat{T}^{\textrm{F}, \vartheta}_k 
\right\Vert_2
&= 
\sqrt{
1 + | \lambda_k |^2 
}
\left| 
\frac{\lambda_k^2 - \mu_k}{\exp(\imath \vartheta) - \mu_k} 
\right|
\leq
\sqrt{
1 + | \lambda_k |^2 
}
\left| 
\frac{\lambda_k^2 - \mu_k}{1 - |\mu_k|} 
\right|
\approx 
\left\Vert
\widetilde{T}^{\textrm{F}}_k 
\right\Vert_2,
\label{eq:LFA_F_T_est}
\\
\left\Vert
\widehat{T}^{\textrm{FCF}, \vartheta}_k 
\right\Vert_2
&= 
\left| \lambda_k \right|^2
\left\Vert
\widehat{T}^{\textrm{F}, \vartheta}_k 
\right\Vert_2
\leq 
\left| \lambda_k \right|^2
\sqrt{
1 + | \lambda_k |^2 
}
\left| 
\frac{\lambda_k^2 - \mu_k}{1 - |\mu_k|} 
\right|
\approx 
\left\Vert
\widetilde{T}^{\textrm{FCF}}_k 
\right\Vert_2.
\label{eq:LFA_FCF_T_est}
\end{align}

It is useful to consider these new LFA estimates in conjunction with the error estimates from \cite{Dobrev2017}. Since the estimates from \cite{Dobrev2017} will also be used in \S \ref{sec:psi_choice}, we reproduce them here. Assuming $|\lambda_k|, |\mu_k| < 1$, and using our current notation, they are \cite[eq. (3.18) \& (3.19)]{Dobrev2017}:
\begin{align}
\left\Vert
\widetilde{T}^{\textrm{F}}_{\Delta,k} 
\right\Vert_2
&\leq
\left| \lambda_k^m - \mu_k \right|
\left(
\frac{1 - \left| \mu_k \right|^{n_t/m}}{1 - \left| \mu_k \right| }
\right)
\quad 
&\textrm{(for F-relaxation),} 
\label{eq:Dobrev_F}
\\
\left\Vert
\widetilde{T}^{\textrm{FCF}}_{\Delta ,k} 
\right\Vert_2
&\leq
\left| \lambda_k \right|^m
\left| \lambda_k^m - \mu_k \right| 
\left(
\frac{1 - \left| \mu_k \right|^{n_t/m-1}}
{ 1 - \left| \mu_k \right| }
\right)
\quad 
&\textrm{(for FCF-relaxation)}.
\label{eq:Dobrev_FCF}
\end{align}
Here $\widetilde{T}_{\Delta,k} \in \mathbb{C}^{n_t/m \times n_t/m}$ are blocks akin to those of \eqref{eq:MGRIT_T_tilde_block}, but describe error propagation on the coarse grid rather than the fine grid; they bound the fine-grid error propagation matrices by $\Vert \widetilde{T}_{k} \Vert_2 \leq \sqrt{m} \Vert \widetilde{T}_{\Delta ,k} \Vert_2$ (see \cite[lemma 4.1]{Hessenthaler2018}).  Additionally, it should be noted that the analysis in \cite{Dobrev2017} was not restricted to circulant $\Phi$ and $\Psi$, but holds for any $\Phi$ and $\Psi$ that are simultaneously unitarily diagonalizable. 

Ignoring simplifications arising from considering only $m = 2$ coarsening in the LFA estimates, we see that \eqref{eq:LFA_F_T_est} and \eqref{eq:LFA_FCF_T_est} share common terms with \eqref{eq:Dobrev_F} and \eqref{eq:Dobrev_FCF}, especially in the limit $n_t \rightarrow \infty$. It is interesting that we can essentially arrive at very similar results to those of \cite{Dobrev2017}, but using a different analysis technique. Given the similarities between the two sets of estimates, and the fact that \eqref{eq:Dobrev_F} and \eqref{eq:Dobrev_FCF} were proven to be tight in \cite{Southworth2018}, it is unsurprising that we have previously observed similar short-term MGRIT convergence between initial-value and time-periodic problems, since \eqref{eq:LFA_F_T_est} and \eqref{eq:LFA_FCF_T_est} apply rigorously for the latter.

\subsection{Discussion}
\label{ssec:analysis_discussion}

Given the new LFA error estimates and those from \cite{Dobrev2017} presented in the previous section, the question is now: What is required of $\Psi$ for effective Parareal/MGRIT convergence?
\begin{enumerate}
\item Clearly in the limiting case that $\Psi = \Phi^m$ the algorithm is exact in one step; however, this is not practically feasible. Nonetheless, observe from the error estimates that convergence of a given spatial mode is dependent on how closely $\mu_k \approx \lambda_k^m$. Thus, in general, for fast convergence, the spectrum of $\Psi$ should approximate that of $\Phi^m$ in some sense.
\item Observe from the denominators in the error estimates that modes for which $|\mu_k| \approx 1$ are potentially damped much slower than those for which $|\mu_k| \ll1$. From the first point, eigenvalues of $\Psi$ should be some approximation to those of $\Phi^m$, and thus, modes for which $|\mu_k| \approx 1$ are associated with $|\lambda_k| \approx 1$. Typically in the context of (dissipative) PDE discretizations, smaller eigenvalues $\lambda_k$ of $\Phi$ are associated with spatially oscillatory modes, $|\theta_k| \approx \pi$, whilst larger eigenvalues are associated with spatially smooth modes, $|\theta_k| \approx 0$. Thus, it can be expected that, in general, convergence of spatially smooth modes is more difficult than for oscillatory modes.
\item Modes with smaller $|\lambda_k|$ are damped much faster by additional relaxation (i.e., FCF over F) compared to those with $|\lambda_k| \approx 1$. Again, this means that it is spatially oscillatory modes that benefit most from additional relaxation.
\end{enumerate}

Given the points above, it is apparent that for fast convergence across all spatial modes the spectrum of $\Phi^m$ should be approximated by that of $\Psi$, and that larger components of the spectrum  (i.e., $|\lambda_k| \approx 1$) should be approximated more accurately. Considering that larger components of the spectrum of $\Phi$ are correlated with smoother spatial components, this equivalently means that the action of $\Phi^m$ on smooth vectors should be more accurately approximated by the action of $\Psi$ on such vectors compared with more oscillatory ones. These conclusions have also been reached in \cite{Southworth2018} using a more general convergence framework.

From this discussion alone, it is not entirely clear why these algorithms tend to perform significantly better for parabolic problems relative to hyperbolic ones. However, from this analysis it is clear that one factor likely aiding in the fast convergence for some schemes is dissipative fine- and coarse-grid time steppers, $\Phi$ and $\Psi$. We remark that discretizations of parabolic PDEs are naturally much more dissipative than those of hyperbolic problems since parabolic PDEs are themselves diffusive, whilst hyperbolic PDEs are advective.

\section{Selection of coarse-grid time stepper $\Psi$} \label{sec:psi_choice}

Based on the conclusions from \S \ref{sec:convergence}  about convergence of Parareal and MGRIT, we now explore selection strategies for the coarse-grid time stepper $\Psi$. We target a $\Psi$ that yields improved convergence relative to simple rediscretization, but whose action is significantly less expensive to compute than the ideal coarse-grid time stepper, $\Phi^m$; the latter criteria is enforced through restrictions on the sparsity of $\Psi$.

Ideally, we seek a $\Psi$ whose eigenvalues minimize an approximation of the two-norm of the error propagator, $
\Vert T \Vert_2 = \max_k 
\big\Vert
\widetilde{T}_{k}
\big\Vert_2$ (with $\widetilde{T}_{k}$ defined in \eqref{eq:MGRIT_T_tilde_block}).
However, this minimax problem is quite difficult to solve, and so we approximate it with the following simpler nonlinear least squares problem:
\begin{align} \label{eq:CG_choice_NLS}
\Psi 
\coloneqq 
\argmin_{\Psi \in \mathbb{R}^{n_x \times n_x}} 
\Vert \bm{p} \Vert^2_2,
\quad
\Vert \bm{p} \Vert^2_2 
\coloneqq 
\sum \limits_{k = -n_x/2}^{n_x/2-1}
\big\Vert
\widetilde{T}_{k,\Delta}
\big\Vert_2^2,
\end{align}
with $\big\Vert \widetilde{T}_{k,\Delta} \big\Vert_2$ being one of the bounds \eqref{eq:Dobrev_F} or \eqref{eq:Dobrev_FCF} from \cite{Dobrev2017}. To arrive at this objective function, we start by minimizing the squared sum of fine-grid iteration matrices \eqref{eq:MGRIT_T_tilde_block} and then bound this by $m \Vert \bm{p} \Vert_2^2$ using $\big\Vert \widetilde{T}_{k}
\big\Vert_2^2 \leq m \big\Vert
\widetilde{T}_{k,\Delta} \big\Vert_2^2$; however, the factor of $m$ is neglected here since it has no effect on the location of minima of $\Vert \bm{p} \Vert_2^2$. We elect to use the bounds from \cite{Dobrev2017} rather than the LFA estimates \eqref{eq:LFA_F_T_est} or \eqref{eq:LFA_FCF_T_est} in this particular case because the LFA estimates are singular for $|\mu_k| \rightarrow 1$, which is likely to be important for the purely advective problem \eqref{eq:LA_PDE}. 

Before discussing further the solution of \eqref{eq:CG_choice_NLS} in \S \ref{ssec:Psi_nonlinear}, an approximate linear problem is explored in \S \ref{ssec:Psi_linear_approx}. Throughout the remainder of this section, we exploit that the eigenvalues of a circulant matrix are simply given by the action of $\sqrt{n_x} {\cal F}$ on its first column (equivalent to the discrete Fourier transform of its first column).

\subsection{Linear least squares approximation for sparse $\Psi$}
\label{ssec:Psi_linear_approx}

Observe from error estimates \eqref{eq:Dobrev_F} and \eqref{eq:Dobrev_FCF} that a first approximation to the solution of \eqref{eq:CG_choice_NLS} might come from minimizing the difference between the spectrums of $\Phi^m$ and $\Psi$. However, given the discussion in \S \ref{ssec:analysis_discussion}, we recognise that it is important to more accurately match some parts of the spectrum of $\Phi^m$ than others. Hence, we propose considering a coarse-grid time stepper defined as the solution of the optimization problem:
\begin{align} \label{eq:CG_choice_min_spectra}
\Psi 
\coloneqq 
\argmin_{\Psi \in \mathbb{R}^{n_x \times n_x}} \Vert \bm{q} \Vert^2_2,
\quad
\Vert \bm{q} \Vert^2_2 
\coloneqq
\sum \limits_{k = -n_x/2}^{n_x/2-1} w_k \left|  \lambda^m_k - \mu_k \right|^2,
\end{align}
in which $w_k \equiv w_k(\lambda_k) > 0$ is a (real) positive weight depending on $\lambda_k$. If $\Psi$ is chosen to be a sparse matrix---corresponding to using explicit time marching on the coarse grid---then \eqref{eq:CG_choice_min_spectra} is a linear least squares problem, of which the solution is now explored.

Let $\hat{\bm{\phi}}^m$, and $\hat{\bm{\psi}} \in \mathbb{R}^{n_x}$
denote the first columns of the matrices $\Phi^m$ and $\Psi$, respectively. Assuming the sparsity pattern of $\Psi$ is given, we let $R_{\Psi} \in \mathbb{R}^{\nu \times n_x}$ be the restriction operator that selects these $\nu$ non-zero entries from $\hat{\bm{\psi}}$. We now define the vector of unknowns 
$\bm{\psi} 
\coloneqq 
R_{\Psi} \hat{\bm{\psi}} 
\in 
\mathbb{R}^{\nu}$
to be the non-zero components of $\hat{\bm{\psi}}$. Finally, let $W \coloneqq \diag \left(w_{-n_x/2}, \dots, w_{n_x/2-1} \right) \in \mathbb{R}^{n_x \times n_x}$ be the weighting matrix.


Using these quantities, \eqref{eq:CG_choice_min_spectra} is written as a linear least squares problem for $\bm{\psi}$:
\begin{align} \label{eq:Psi_selection_LLS_problem}
\bm{\psi} 
\coloneqq 
\argmin \limits_{\bm{\psi} \in \mathbb{R}^{\nu}}  
\left\Vert 
W^{1/2} 
{\cal F}
\left(
\hat{\bm{\phi}}^m - R_{\Psi}^\top \bm{\psi} 
\right)
\right\Vert_2^2.
\end{align}
Forming and (symbolically) solving the normal equations for this problem we find
\begin{align} \label{eq:linlsq_gen_sol}
\bm{\psi} 
= 
\left( W^{1/2} {\cal F} R_{\Psi}^\top \right)^+ W^{1/2} {\cal F} \hat{\bm{\phi}}^m,
\quad
\hat{\bm{\psi}} 
= 
R_{\Psi}^\top \bm{\psi},
\end{align}
in which $X^+ \equiv \left( X^* X \right)^{-1} X^*$ denotes the pseudoinverse of matrix $X$. One special case to consider is for weights $w_k  = 1$, which corresponds to minimizing the difference between the spectra of $\Phi^m$ and $\Psi$ in the two-norm. For this special case, the solution can be simplified and has a straightforward interpretation: $\Psi$ is given by truncating $\Phi^m$ in the sparsity pattern of $\Psi$. 

Note that, in general, the solution \eqref{eq:linlsq_gen_sol} is not guaranteed to be real; for results reported in \S \ref{ssec:Psi_results}, the imaginary components of the solution are simply discarded since they are found to be many orders of magnitude smaller than the real components.


\subsection{Nonlinear least squares solution}
\label{ssec:Psi_nonlinear}

We now discuss the solution of the nonlinear least squares problem \eqref{eq:CG_choice_NLS}. Note that we are not necessarily advocating that the approach described here for selecting $\Psi$ is practically feasible: our primary concern is in trying to better understand what sorts of convergence rates \textit{are} possible.

Problem \eqref{eq:CG_choice_NLS} is solved using the nonlinear least squares routine \texttt{lsqnonlin} from MATLAB's optimization toolbox. To do so, \texttt{lsqnonlin} needs to be able to evaluate the residual vector $\bm{p}$ in \eqref{eq:CG_choice_NLS}, which requires the eigenvalues $\lambda_k$ and $\mu_k$. The particular form of these eigenvalues depends on the form of $\Phi$ and $\Psi$, respectively. In the most general case considered here, they are each the product of the inverse of a sparse matrix and a sparse matrix:
\begin{align} \label{eq:CG_choice_general_Phi_and_Psi}
\Phi \coloneqq \Phi_{\rm I}^{-1} \Phi_{\rm E},
\quad
\Psi \coloneqq \Psi_{\rm I}^{-1} \Psi_{\rm E}.
\end{align}
The sparse matrices $\Psi_{\rm I}$, and $\Psi_{\rm E}$ (I and E respectively denoting implicit and explicit) are to be determined, but their sparsity patterns are assumed to be specified a priori. We noted in \S \ref{ssec:model_problem} that $\Phi$ takes the form of a rational function in the spatial discretization. This behaviour is realized in \eqref{eq:CG_choice_general_Phi_and_Psi} through $\Phi_{\rm I}$ and $\Phi_{\rm E}$ playing the roles of the denominator and numerator, respectively, of this rational function, and thus being polynomials in the spatial discretization. The polynomials may be determined by appealing to the stability function of the given RK method, as discussed in \cite{Hessenthaler2018}.

Analogously to \S \ref{ssec:Psi_linear_approx}, we define $\hat{\bm{\phi}}_{\rm I}$, $\hat{\bm{\phi}}_{\rm E}$, $\hat{\bm{\psi}}_{\rm I}$, and $\hat{\bm{\psi}}_{\rm E}$ to be the first columns of $\Phi_{\rm I}$, $\Phi_{\rm E}$, $\Psi_{\rm I}$, and $\Psi_{\rm E}$, respectively. Using  the sparsity patterns of $\Psi_{\rm I}$ and $\Psi_{\rm E}$, we define the restriction operators $R_{\Psi_{\rm I}} \in \mathbb{R}^{\nu_{\rm I} \times n_x}$, and $R_{\Psi_{\rm E}} \in \mathbb{R}^{\nu_{\rm E} \times n_x}$ that select the non-zero elements from $\Psi_{\rm I}$ and $\Psi_{\rm E}$, respectively. Using these restriction operators, the vectors of unknowns are defined as 
$\bm{\psi}_{\rm I} 
\coloneqq 
R_{\Psi_{\rm I}} \hat{\bm{\psi}}_{\rm I} 
\in 
\mathbb{R}^{\nu_{\rm I}}, $ and $
\bm{\psi}_{\rm E} 
\coloneqq 
R_{\Psi_{\rm E}} \hat{\bm{\psi}}_{\rm E} 
\in 
\mathbb{R}^{\nu_{\rm E}}$.
The $k$th eigenvalues of the time steppers defined in \eqref{eq:CG_choice_general_Phi_and_Psi} are thus given by
\begin{align} 
\lambda_{k} = 
\frac{ 
\big( {\cal F} \hat{\bm{\phi}}_{\rm E} 
\big)_k 
}
{ 
\big( {\cal F} \hat{\bm{\phi}}_{\rm I} 
\big)_k 
},
\quad
\mu_{k} = 
\frac{ 
\big( 
{\cal F} R_{\Psi_{\rm E}}^\top \bm{\psi}_{\rm E} 
\big)_k 
}
{ 
\big( 
{\cal F} R_{\Psi_{\rm I}}^\top \bm{\psi}_{\rm I} 
\big)_k 
}.
\end{align}
For $\Phi$ and $\Psi$ in the form of \eqref{eq:CG_choice_general_Phi_and_Psi}, $\bm{p}$ in \eqref{eq:CG_choice_NLS} can be computed using these expressions.

\subsection{Results}
\label{ssec:Psi_results}

Numerical results obtained from applying the coarse-grid selection strategies of the previous sections to the discretizations described in \S \ref{ssec:model_problem} are now presented. In all experiments, the initial guess at the space-time solution is taken to be uniformly random for all times $t >0$. The experimental convergence metric used is the number of iterations required to reduce the two-norm of the space-time residual vector from its initial value by 10 orders of magnitude. We search for components $\Psi_{\rm I}$ and $\Psi_{\rm E}$ of the coarse-grid time stepper \eqref{eq:CG_choice_general_Phi_and_Psi} having the same sparsity patterns as their fine-grid counterparts, $\Phi_{\rm I}$ and $\Phi_{\rm E}$, respectively. Thus, for the two explicit schemes we target sparse coarse-grid time steppers ($\Psi_{\rm I} \equiv I$) and so we use both the weighted linear least squares approach of \S \ref{ssec:Psi_linear_approx} and the nonlinear strategy of \S \ref{ssec:Psi_nonlinear}. \renewcommand*{\arraystretch}{1.2}
\begin{table}[tbhp] \label{tb:explicit}
{\footnotesize
\caption{Iteration counts for 2nd- and 3rd-order explicit schemes. {\rm weighted lin.} and {\rm nonlin.} denote results with $\Psi$ determined with the weighted linear and nonlinear least squares strategies, respectively.}
\begin{center}
\begin{tabular}{|c|c|c|c|c|c|} 
\cline{2-5}
\multicolumn{1}{c|}{} & \multicolumn{2}{|c|}{2nd-order}  & \multicolumn{2}{c|}{3rd-order}  \\ \hline
$n_x \times n_t$ 	& weighted lin. & nonlin. & weighted lin. & nonlin. \\ \hline
$2^6 \times 2^6$ 		& 11 & 9  & 13 & 11 \\ \hline
$2^8 \times 2^8$ 		& 11 & 9 & 15 & 11 \\ \hline
$2^{10} \times 2^{10}$ 	& 11 & 9 & 15 & 11 \\ \hline
$2^{12} \times 2^{12}$ 	& 11 & 9 & 15 & 11 \\ \hline
\end{tabular}
\end{center}
}
\end{table}
\renewcommand*{\arraystretch}{1.2}
\begin{table}[tbhp] \label{tb:implicit}
{\footnotesize
\caption{Iteration counts for 3rd-order implicit scheme. {\rm redisc.} and {\rm nonlin.} denote results with $\Psi$ given by rediscretization and the nonlinear least squares strategy, respectively. {\rm DNC} denotes a solve that did not converge to the desired tolerance within 50 iterations.}
\begin{center}
\begin{tabular}{|c|c|c|c|c|c|c|c|c|c|}
\cline{2-9} 
\multicolumn{1}{c|}{} 
& \multicolumn{4}{c}{F-relaxation} & \multicolumn{4}{|c|}{FCF-relaxation} \\
\cline{2-9}  
\multicolumn{1}{c|}{}  
& \multicolumn{2}{|c|}{$m=2$} & \multicolumn{2}{|c|}{$m=4$} & \multicolumn{2}{|c|}{$m=2$} & \multicolumn{2}{|c|}{$m=4$}  \\ \hline
$n_x \times n_t$ & redisc. & nonlin. & redisc. & nonlin. & redisc. & nonlin. & redisc. & nonlin.\\ \hline
$2^6 \times 2^6$ 		& 20  & 5 & 16     & 4 & 14 &  5 & 9       & 5  \\ \hline
$2^8 \times 2^8$ 		& 31  & 5 & DNC & 6  & 24 & 5 & 32     & 5 \\ \hline
$2^{10} \times 2^{10}$ 	& 35 & 5 & DNC & 6  & 31 &  5 & DNC & 5 \\ \hline
$2^{12} \times 2^{12}$ 	& 36 & 5 & DNC & 6  & 32 & 5 & DNC & 5 \\ \hline
\end{tabular}
\end{center}
}
\end{table}

For the explicit schemes, we use FCF relaxation with $m = 2$ coarsening. The 2nd-order scheme is run at a CFL number of $a \tfrac{\Delta t}{\Delta x} = 0.4$ and uses linear least squares weights $w_k = |\lambda_k|^{40}$, whilst the 3rd-order scheme uses a CFL number of 1.4 and weights $w_k = |\lambda_k|^{20}$. These choices of weights are somewhat arbitrary, but are motivated by the fact that there should be relatively close agreement between the larger elements of the spectra of $\Phi^m$ and $\Psi$ (\S \ref{ssec:analysis_discussion}), and because experimentally they yielded scalable solvers (at least up to the maximum grid size considered here). The results of these experiments are given in Table \ref{tb:explicit}. All resulting solvers have quite fast convergence rates. For both the 2nd- and 3rd-order cases, the weighted linear least squares strategy does not yield as fast a solver as the nonlinear one (as is to be expected), but nonetheless it provides an excellent approximation.

Note that iteration counts for solvers using rediscretization coarse-grid operators are not included in Table \ref{tb:explicit}. Explicit fine-grid discretizations that use such coarse-grid operators are imposed with an overly restrictive fine-grid CFL condition to ensure a stable coarse-grid discretization. This restriction typically makes parallel time integration for explicit discretizations intractable, since the discretizations themselves are only considered effective when run at a large fraction of their fine-grid CFL limit. A promising feature of the strategies presented here is that they do not use rediscretization, and so they do not necessarily suffer from this instability. In particular, this is true for the explicit results shown in Table \ref{tb:explicit}, where effective convergence rates are achieved using fine-grid CFL numbers that are significant fractions of their respective fine-grid CFL limits (80\% and $86\%$ for the 2nd- and 3rd-order schemes, respectively). Consequently, the fine-grid CFL numbers significantly exceed those needed to obtain stable rediscretizations on their respective coarse grids (meaning solvers with rediscretization coarse-grid operators here would be divergent). These results represent a novel contribution to the solution of explicit time discretizations with MGRIT.

We run experiments for the 3rd-order implicit scheme using both rediscretization and nonlinear least squares coarse-grid operators. We use both $m = 2$ and $m = 4$ coarsening, both of which run a fine-grid CFL number of unity. The iteration counts for these experiments are given in Table \ref{tb:implicit}. Rediscretization yields a scalable solver for $m = 2$, but convergence is quite slow. For the $m=4$ case, rediscretization is essentially divergent: convergence is only really achieved for iteration counts approximately equal to those for which the algorithms reproduce the sequential time marching solution \cite{Falgout2014}. However, the nonlinear coarse-grid operator yields a scalable solver in all cases with rapid convergence, and with only very small differences between $m = 2$ and $m = 4$ coarsening, and between F- and FCF-relaxation. 

We conclude this section by remarking that the rapid and scalable convergence achieved here for both (relatively non-dissipative) high-order implicit and explicit discretizations of the advection equation using Parareal and MGRIT represent a novel contribution to the field. Much of the existing literature \cite{DeSterck2018_mgrit_hyp, HowseThesis2017, Ruprecht2018} has focused on dissipative 1st-order discretizations for which large iteration counts are observed and results are not scalable with problem size. Note that in \cite{Dobrev2017}, high-order implicit discretizations were considered, but yielded slower convergence rates than presented here and were not scalable without the use of artificial dissipation.

\section{Conclusions}\label{sec:conclusions}
To better understand the convergence behaviour of Parareal and two-level MGRIT algorithms, we have reviewed some existing convergence theory and developed  new theory using Fourier analysis. This new convergence theory is the first contribution of this paper. Motivated by this convergence theory, our second contribution is developing selection strategies for coarse-grid operators, which can, in principle, be used to create efficient algorithms. This is significant because it demonstrates that there exist coarse-grid operators that lead to efficient and scalable Parareal and MGRIT algorithms for problems where existing coarse-grid operators (i.e., rediscretization) do not. Our model test problem was the linear advection equation, which is well known to cause significant difficulties for these algorithms, even in cases of simple 1st-order discretizations. We demonstrate there exist coarse-grid operators for which effective convergence rates are achievable for this problem. Specifically, for the 2nd- and 3rd-order explicit discretizations considered, we identified scalable solvers with relatively fast convergence rates. Most notably, this was done for problems run at significant fractions of their fine-grid CFL limits, and thus corresponds to cases in which rediscretization coarse-grid operators are divergent due to their violating of a coarse-grid CFL condition. For a 3rd-order implicit discretization, rapid and scalable convergence was attained for cases in which rediscretization coarse-grid operators were either very slow to converge or divergent. The proposed coarse-grid selection strategies thus lead to significant improvements over existing techniques, which we hope will lead to achieving greater parallel speed-up for the advection equation in the future.

Future work includes developing a better understanding of cases for which the proposed strategies fail to yield improved convergence over rediscretization, as is the case for 1st-order discretizations of the model problem, for example. Developing better weights for the proposed least squares strategy, and generalizing these ideas to discretizations with different boundary conditions will also be considered. 

\bibliographystyle{siamplain}
\bibliography{references}

\end{document}